\documentclass[11pt,oneside]{article}
\usepackage{amsfonts,amssymb}
\usepackage[centertags]{amsmath}
\usepackage[italian,english]{babel}
\usepackage{theorem}
\usepackage[dvips]{graphicx}
\newcommand{\R}{\ensuremath{\mathbb R}}
\newcommand{\C}{\ensuremath{\mathbb C}}
\renewcommand{\Re}{\mathop{\mathrm{Re}}}
\renewcommand{\Im}{\mathop{\mathrm{Im}}}
\newcommand{\rank}{\mathop{\mathrm{rank}}}
\newcommand{\Ad}{\mathop{\mathrm{Ad}}}
\newcommand{\End}{\mathop{\mathrm{End}}}

\renewcommand{\setminus}{\smallsetminus}
\newcommand{\tr}{\mathop{\mathrm{tr}}}
\renewcommand{\ge}{\geqslant}
\renewcommand{\le}{\leqslant}

\renewcommand{\j}{{\ol j}}
\def\1{{\ol1}}\def\2{{\ol2}}\def\3{{\ol3}}
\newcommand{\ba}{\begin{array}}
\newcommand{\ea}{\end{array}}
\newcommand{\vs}{\vspace{5pt}}
\newcommand{\y}{\\[4pt]}
\newcommand{\yy}{\\[6pt]}
\newcommand{\ol}{\overline}
\newcommand{\ti}{\tilde}
\newcommand{\vph}{\vphantom{\int^{1^1}}}
\newcommand{\q}{\quad}\newcommand{\qq}{\qquad}
\newcommand{\cC}{\mathcal{C}}
\newcommand{\cG}{\mathcal{G}}
\newcommand{\g}{\mathfrak{g}}
\newcommand{\n}{\noindent}
\newcommand{\we}{\wedge}
\newcommand{\bs}{\backslash}
\newcommand{\na}{\nabla}
\newcommand{\al}{\alpha}
\newcommand{\om}{\omega}
\newcommand{\ga}{\gamma}
\newcommand{\de}{\delta}
\newcommand{\Ga}{\Gamma}

\newcommand{\Om}{\Omega}
\newcommand{\pd}{\partial}
\newcommand{\opd}{\ol\pd}
\newcommand{\la}{\lambda}
\newcommand{\La}{\Lambda}
\newcommand{\si}{\sigma}
\renewcommand{\th}{\theta}
\newcommand{\B}{\mathbf{B}}
\newcommand{\DD}{\mathbb{D}}
\newcommand{\I}{\mathbf{I}}
\newcommand{\cI}{\mathcal{I}}
\newcommand{\M}{\mathbf{M}}

\renewcommand{\P}{\mathbf{P}}
\newcommand{\Q}{\mathbf{Q}}
\newcommand{\X}{\mathbf{X}}
\newcommand{\Y}{\mathbf{Y}}
\newcommand{\Z}{\mathbf{Z}}

\newcommand{\hJ}{\hat J}
\newcommand{\XoX}{\X\ol\X}
\newcommand{\alt}{\raise1pt\hbox{$\bigwedge$}}
\newcommand{\cir}{\hbox{\,\footnotesize$\circ$}\,}
\newcommand{\Minors}{\mathrm{Minors}}
\newcommand{\fs}[2]{\hbox{\large$\frac{#1}{#2}$}}
\newcommand{\ft}[2]{\hbox{$\frac{#1}{#2}$}}
\renewcommand{\sb}[1]{{}\!_{#1}}
\renewcommand{\a}[1]{\al^{#1}}
\renewcommand{\aa}[2]{\a{#1#2}}
\newcommand{\ao}[2]{\a{#1\ol#2}}
\newcommand{\oa}[2]{\a{\ol#1#2}}
\newcommand{\w}[1]{\om^{#1}}
\newcommand{\ww}[2]{\w{#1#2}}
\newcommand{\wo}[2]{\w{#1\ol#2}}
\newcommand{\ow}[2]{\w{\ol#1#2}}
\newcommand{\xo}[2]{x_{#1\ol#2}}

\newcommand{\sq}{\kern.5pt\square}
\newtheorem{lem}{Lemma}[section]

\newtheorem{teo}[lem]{Theorem}
\newtheorem{cor}[lem]{Corollary}
\newtheorem{pro}[lem]{Proposition}
{\theorembodyfont{\rmfamily} 
 \newtheorem{oss}[lem]{Remark}
 \newtheorem{defi}[lem]{Definition}
 \newtheorem{exa}[lem]{Example} }
\newcommand{\hfr}{\hbox{$\frown$}}
\newcommand{\hsm}{\hbox{$\smile$}}
\newcommand{\qsm}{\raise1.5pt\hfr\kern-11pt\lower1.5pt\hsm}
\newcommand{\EndDim}{\ensuremath{\nopagebreak\hfill\qsm}\medbreak}

0
\newenvironment{D}[1][]{{\nopagebreak\n\em Proof{#1}: }}{\EndDim}

\title{Families of strong KT structures in six dimensions}
\author{Anna Fino, Maurizio Parton and Simon Salamon}
\date{}
\begin{document}

\maketitle

\n{\small\emph{Abstract.} This paper classifies Hermitian structures on
6-dimensional nilmanifolds $M=\Ga\bs G$ for which the fundamental 2-form is
$\pd\opd$-closed, a condition that is shown to depend only on the underlying
complex structure $J$ of $M$. The space of such $J$ is described when $G$ is
the complex Heisenberg group, and explicit solutions are obtained from a lima\c
con-shaped curve in the complex plane. Related theory is used to provide
examples of various types of Ricci-flat structures.}\vs


\section*{Introduction}

Let $(M,J,g)$ be a Hermitian manifold. There is a 1-parameter family of
canonical Hermitian connections on $M$ which can be distinguished by properties
of the torsion tensor $T$ \cite{GauHCD,YanDGC}. In particular, there is a
unique connection $\na$ satisfying $\na g=0$, $\na J=0$ for which $g(X,T(Y,Z))$
is totally skew-symmetric. The resulting 3-form can then be identified with
$Jd\Om$, where $\Om$ is the fundamental 2-form defined by \eqref{gom}. This
connection was used by Bismut \cite{BisLIT} to prove a local index formula for
the Dolbeault operator when the manifold is non-K\" ahler, and is the subject
of a number of other interesting results \cite{GIPDGC,IvPVTS}. The properties
of such a connection give rise to what is loosely called `K\"ahler with torsion
geometry', and if $Jd\Om$ is closed but non-zero then $g$ is called a
\emph{strong KT metric}. Such metrics have applications in type II string
theory and in 2-dimensional supersymmetric $\sigma$-models
\cite{GHRTMN,HoPFRG,StrSWT}.

In four real dimensions, a metric satisfying the strong KT condition is
`standard' in the terminology of Gauduchon \cite{Gau1FT}. One can be found in
the conformal class of any given Hermitian metric on a compact manifold. But
the theory is very different in higher dimensions. Even-dimensional compact
Lie groups provide a natural class of strong KT structures \cite{SSTCSP}. In
this case, one may choose $J$ to be a left-invariant complex structure and $g$
to be a compatible bi-invariant metric. Then $\na$ is the flat connection with
skew-symmetric torsion $g(X,[Y,Z])$ corresponding to an invariant 3-form on
the Lie algebra. It is therefore natural to investigate the situation with
regard to other groups.

In this paper, we study KT geometry on 6-dimensional nilmanifolds in which $J$
and $g$ arise from corresponding left-invariant tensors. If $G$ is a
simply-connected nilpotent Lie group, and if the structure equations of its Lie
algebra are rational, then there exists a discrete subgroup $\Ga$ of $G$ for
which $M=\Ga\bs G$ is compact \cite{MalCHS,NomCCH}. Any left-invariant complex
structure on $G$ will pass to a complex structure $J$ on $M$ but, unless $G$ is
abelian, the $\pd\opd$-lemma is not valid for $J$ and in particular there is
no compatible K\"ahler metric \cite{BeGKSS,DGMRHT,HasMMN,NeTDHT}. As we explain
in \S6, there may or may not be invariant pseudo-K\"ahler metrics on $(M,J)$.

Eighteen of the thirty-four classes of real 6-dimensional nilpotent Lie
algebras $\g$ admit a complex structure. Exactly four of these classes, all of
them 2-step with $b_1\ge4$ and including the case in which $\g$ underlies the
complex Heisenberg algebra, give rise to strong KT metrics. Given that compact
nilmanifolds with a strong KT structure exist, it is perhaps surprising that
there are so few classes. The classification over $\R$ is accomplished in \S3,
after an analysis of the relevant structure equations over $\C$ in \S\S1,2. A
matrix formalism for describing (1,1)-forms is introduced in an attempt to
make the calculations of this paper rather more enlightening. A striking
feature of our classification is that the existence of a strong KT structure
depends only on the complex structure of $\g$, and this poses the question of
understanding the solutions as a subset of an appropriate moduli space of
complex structures.

With this aim, we proceed to a detailed study of the strong KT equations when
$G$ is the complex Heisenberg group and $M=\Ga\bs G$ is the Iwasawa
manifold. It is easy to check that none of the standard complex structures
\cite{AGSAHG} on $G$ are strong KT, so we were intrigued to discover which ones
are. According to the third author's joint paper with Ketsetzis \cite{KeSCSI},
essential features of an invariant complex strcture $J$ on $M$ depend on
$\XoX$, where $\X$ is a $2\times2$ matrix representing the induced action of
$J$ on $M/T^2\cong T^4$. In \S4, we prove that the strong KT condition
constrains the eigenvalues of $\XoX$ to be complex conjugates lying on a curve
in the complex plane. We interpret this result in terms of the action of the
automorphism group of $\g$ in \S5, and this leads to an explicit description of
the solution space. An analogous study can probably be carried out when
$G=H_3\times H_3$ is the product of real Heisenberg groups, using methods from
\cite{GPPKOD}.

A Hermitian manifold is called conformally balanced if the Lee 1-form $\theta$
(the `trace' of $d\Om$) is exact. The study of such structures in connection
with the connection $\na$ is motivated by work of \cite{PapKHG}, though there
are less subtleties in our context in which $\th$ is exact if and only if it is
zero. Having explained in \S1 that the vanishing of $\th$ is complementary to
the SKT condition, we observe in the final section that $\th$ is also the
obstruction to the holonomy of $\na$ reducing to $SU(n)$. We list some
6-dimensional Lie algebras giving rise to nilmanifolds admitting such a
reduction, and others admitting a psuedo-Riemannian metric with zero Ricci
tensor.\vs

\n{\small Acknowledgment. The authors are members of EDGE, Research Training
Network HPRN-CT-200-00101, supported by the European Human Potential
Programme.}

\section{Complex structure equations}

Let $(M,J,g)$ be a Hermitian manifold of real dimension $2n$. We shall regard
the complex structure $J$ as the primary object, so the Riemannian metric $g$
is chosen to render $J$ orthogonal. The fundamental 2-form $\Om$ is then
defined by \begin{equation}\label{gom} g(X,Y)=\Om(X,JY)\end{equation} and has
type (1,1) relative to $J$. The Hermitian structure is K\"ahler if and only if
$d\Om=0$, which is equivalent to the vanishing of $\pd\Om=(d\Om)^{2,1}$.

Somewhat unconventionally, we set \[\sq=\ft12idJd.\] This operator acts as
$\pd\opd$ on forms of type $(p,p)$, which it maps to forms of type
$(p+1,p+1)$. We shall only be concerned with the case $p=1$.

\begin{defi}\label{strong} We shall say that the Hermitian manifold $(M,J,g)$ 
is `strong KT' or more briefly `SKT' if $\sq\Om=0$ but $d\Om\ne0$.\end{defi}

\n Observe that our definition of SKT excludes the K\"ahler case.\vs

We wish to combine the notion of SKT with that of an invariant Hermitian
structure on a nilmanifold. First recall the definition of nilpotency for a
Lie algebra $\g$. The descending central series of $\g$ is the chain of ideals
defined inductively by $\g^0=\g$ and $\g^i=[\g^{i-1},\g]$ for $i\ge1$. Then
$\g$ is nilpotent if $\g^s=0$ for some $s$. If, in addition, $\g^{s-1}\ne0$
then $\g$ is said to be $s$-step.

Let $\g$ be a real $2n$-dimensional nilpotent Lie algebra. Assigning an almost
complex structure $J:\g\to\g$ is equivalent to choosing an $n$-dimensional
subspace $\La$ of $\g_c^*$ such that $\La\cap\ol\La=\{0\}$. For the purpose of
this paper, we shall call such a subspace of $\g_c^*$ `maximally complex'.
The endomorphism $J$ extends uniquely to a left-invariant almost complex
structure (also denoted by $J$) on any Lie group $G$ with Lie algebra $\g$. The
subspace $\La$ generates the space of (1,0)-forms relative to $J$, and this is
a complex structure if and only if $\cI(\La)$ is a differential ideal. (We use
$\cI(S)$ to denote the ideal generated by a subset $S$ of the exterior algebra
$\alt^*\g$.)

By a `nilmanifold with an invariant complex structure' we mean an
even-dimensional nilmanifold $\Ga\bs G$ endowed with a complex structure $J$
arising from $\g$. It is important to note that $(G,J)$ will not in general be
a complex Lie group. 

\begin{teo}\label{SKT} Let $M=\Ga\bs G$ be a 6-dimensional nilmanifold with an
invariant complex structure $J$. Then the SKT condition is satisfied by either
all invariant Hermitian metrics $g$ or by none. Indeed, it is satisfied if and
only if $J$ has a basis $(\a i)$ of $(1,0)$-forms such that
\begin{equation}\label{SKT1}\left\{\ba{l}d\a1=0\\d\a2=0\\d\a3=A\oa12+
B\oa22+C\ao11+D\ao12+ E\aa12\ea\right.\end{equation} where $A,B,C,D,E$ are
complex numbers such that \begin{equation}\label{SKT2}
|A|^2+|D|^2+|E|^2+2\Re(\ol BC)=0.\end{equation}\end{teo}

\n We indicate $\a i\we\a\j$ by $\ao ij$ (or $-\aa\j i$), and use similar 
notation for forms of arbitrary degree. Thus, the symbol $\al$ stands more for 
the choice of basis than for an individual element.

The third equation in \eqref{SKT1} means that \[d\a3\in\alt^2\langle\a1,
\a{\ol1},\a2,\a{\ol2}\rangle.\] and the resulting complex structure $J$ is of
`nilpotent' type in the language of \cite{CFGCNN}. Observe that the system
\eqref{SKT1} automatically satisfies $d^2=0$ and therefore defines a Lie
algebra irrespective of the values of $A,B,C,D,E$. The resulting isomorphism
classes are listed in \S3, where we distinguish those compatible with
\eqref{SKT2}.

We shall divide the proof of Theorem~\ref{SKT} into two parts. The second part
is devoted to a mainly computational derivation of the structure equations
\eqref{SKT1}, and is relegated to the next section. However, it is instructive
to begin by assuming \eqref{SKT1} and deducing \eqref{SKT2} from it. This we do
immediately, and in passing we shall see that the choice of metric is
irrelevant.

Let $\Om$ be the fundamental 2-form of some $J$-Hermitian metric, and set
\begin{equation}\label{Omega}\Om=\sum_{i,j=1}^3\xo ij\,\aa i\j,\end{equation}
where $\xo ij\in\C$ are constant coefficients with $\ol{\xo ji}=-\xo ij$. The
positive definiteness of $g$ implies that the restriction of $\Om$ to any
complex line is non-zero.  Equivalently, \begin{equation}\label{restriction}
\Om(V,J\ol V)=g(V,\ol V)>0\end{equation} for any vector $V\ne0$ in the
complexified tangent space.

Given that $Jd\a3$ differs from $d\a3$ by changing the sign of $E$, an easy
calculation gives \begin{equation}\label{ABC} Jd\a3\we d\a\3= (|A|^2+B\ol
C+C\ol B+|D|^2+|E|^2)\a{1\12\2}.\end{equation} The vanishing of this 4-form is
precisely \eqref{SKT2}, which can now be deduced from \eqref{SKT1} via

\begin{lem}\label{vanish} Given \eqref{SKT1}, $(M,J,g)$ is SKT for
any invariant Hermitian metric $g$ if and only if $Jd\a3\we d\a\3=0$.\end{lem}

\begin{D} If $(V_i)$ is a basis of (1,0) vectors dual to $(\a i)$ then
\eqref{restriction} implies that $\xo33=\Om(V_3,J\ol V_3)>0$. On the other
hand, all the terms in \eqref{Omega}, with the exception of $\xo33\ao33$, are
eliminated by two differentiations. Thus, the SKT condition is satisfied for
any compatible metric iff $\sq\ao33=0$.

Let $\Psi=Jd\a3\we d\a\3$, and observe that \[\ba{rcl}2\sq\ao33\>=\>
dJd(i\ao33)&=& idJ(d\a3\we \a\3-\a3\we d\a\3)\y &=&d(Jd\a3\we \a\3+\a3\we
Jd\a\3)\y &=& Jd\a3\we d\a\3+d\a3\we Jd\a\3\y &=&\Psi+J\Psi\ea\] since $J^2=1$
on 2-forms. But $\Psi$ is a form of type (2,2) relative to $J$ (or any other
almost complex structure on the real 4-dimensional space underlying
$\langle\a1,\a2\rangle$), so $J\Psi=\Psi$.\end{D}

The Lee form of a Hermitian manifold $(M,J,g)$ real dimension $2n$ is
the 1-form \[\th=J*d*\Om=-Jd^*\Om\] where $d^*$ is the formal adjoint
of $d$ with respect to $g$. The formula $*\Om=\Om^{n-1}/(n-1)!$
implies that $d(\Om^{n-1}) =\th\we\Om^{n-1}$. Equivalently,
\[\xi=d\Om-\fs1{n-1}\th\we \Om\] satisfies $\xi\we\Om^{n-2}=0$, and is
therefore a \emph{primitive} form. Under a conformal change $\ti g=e^{2f}g$,
the Lee form transforms as \[\ti\th=\th+ 2(n-1)df.\] Almost Hermitian manifolds
with $\th=0$ have in the past been called \emph{semi-K\"ahler} or
\emph{cosymplectic}, though a Hermitian structure is also called
\emph{balanced} if $\th=0$. In this case, we are therefore talking about
Hermitian manifolds of Gray-Hervella class $\mathcal{W}_3$ \cite{GrHSCA}.

The Hermitian structure is \emph{conformally balanced} if $\th$ is exact, for
in that case $f$ can be chosen so that $\ti\th=0$. However, in the invariant
setting, $\th$ is exact if and only if $\th=0$. Since $\langle d^*\Om,\si
\rangle=\langle\Om,d\si\rangle$ for all 2-forms $\si$, the vanishing of $\th$
is equivalent to $\Om$ being orthogonal to the image of $d$ in $\alt^2\g^*$, a
fact exploited in the study \cite{AGSAHG}.

In real dimension 4, the SKT condition is equivalent to $d^*\th=0$. On the
other hand,

\begin{pro}\label{incomp} A Hermitian manifold $(M,J,g)$ of real dimension 
$2n\ge6$ can only be SKT if $\th\ne0$.\end{pro}

\begin{D} Suppose that $\th=0$. Since $d^*\Om=*\,d*\Om$ and $*\Om$ is 
proportional to $\Om^{n-1}$, we may conclude that the $(n-1)$ form
$\Phi=\Om^{n-2}\we\pd\Om$ vanishes. This implies that $\pd\Om$ is a
\emph{primitive} (2,1) form, and \[0=\opd\Phi=
(n-2)\Om^{n-3}\we\opd\Om\we\pd\Om+\Om^{n-2}\we\opd\pd\Om.\]
Primitivity implies that $\Om^{n-3}\we\ol{\pd\Om}\we\pd\Om$ is
proportional to $\|\pd\Om\|^2$, so $\opd\pd\Om=0$ now implies that
$\pd\Om=0$ and $M$ is K\"ahler. But this is excluded in
Definition~\ref{strong}.\end{D}

It is amusing to view this result in the light of Theorem~\ref{SKT}. Suppose
that $n=3$ and that \eqref{SKT1} holds. Referring to \eqref{Omega}, we may set
\[\xo12=z,\q\xo21=-\ol z,\q\xo11=ix,\q\xo22=iy\] with $x,y>0$ and $xy>|z|^2$
to reflect positivity. The condition that $\Om$ be orthogonal to $d\a3$
becomes \begin{equation}\label{bal}Az+i(Bx-Cy)+D\ol z=0.\end{equation}
Proposition~\ref{incomp} implies that this is incompatible with the inequality
\begin{equation}\label{ineq}|A|^2+|D|^2+ 2\Re(\ol BC)\le0\end{equation} from
\eqref{SKT2} unless all the coefficients vanish.

The incompatibility between \eqref{bal} and \eqref{ineq} is clear if $\Om$
assumes the standard form
\begin{equation}\label{Om}\Om_0=\ft12i(\wo11+\wo22+\wo33).\end{equation} For
then $z=0$, $x=y$; thus $B=C$ and $\Re(B\ol C)=|B|^2$. The general case is far
less obvious, but follows by setting $B=y=1$ (which is no real restriction)
and applying

\begin{oss} Let $A,C,D,z\in\C$ be such that $x=C+i(Az+D\ol z)$ is real.
Then \[|A|^2+|D|^2+2\Re C\le0\q\Rightarrow\q x\le|z|^2.\] To verify this, set
$F=i(Az+D\ol z)$ so that \[x=\Re C+\Re F\le\Re C+|F|\le-\ft12(|A|^2+|D|^2)
+(|A|+|D|)|z|.\] If $x>|z|^2$ then \[2|z|^2-2(|A|+|D|)|z|+|A|^2+|D|^2<0\] which
(as a quadratic in $|z|$ with non-positive discriminant) is impossible.
\end{oss}

\section{Reducing the coefficients}

This section is devoted to completing the proof of Theorem~\ref{SKT} by
arriving at \eqref{SKT1}. Our starting point is \cite[Theorem~1.3]{SalCSN},
which we mildly re-state as

\begin{teo}\label{teosal} A maximally complex subspace $\La$ of $\g^*_c$ 
is the $(1,0)$-space of a complex structure if and only if $\La$ has a basis
$(\a i)$ for which $d\a1=0$ and $d\a i\in\cI(\{\a1,\ldots,\a{i-1}\})$ for
$i\ge2$.\end{teo}

Recall that $\cI$ stands for `ideal'. In dimension 6 this result provides the
generic structure equations \begin{equation}\label{generic}\ba{l} d\a1=0,\y
d\a2=a_1\aa12+a_2\aa13+a_3\ao11+a_4\ao12+a_5\ao13,\y
d\a3=b_1\aa12+b_2\aa13+b_3\ao11+b_4\ao12+b_5\ao13+
c_1\aa23+c_2\ao21+c_3\ao22+c_4\ao23\ea\end{equation} where
$a_i,b_j,c_k\in\C$. We therefore need to establish the vanishing of nine
coefficients, namely $a_*$ and $b_2,b_5,c_1,c_4$. (The new names $C,D,-A,-B$
for $b_3,b_4,c_2,c_3$ are carefully chosen to simplify formulae in subsequent
sections.) We follow a type of decision tree in order to eliminate the
coefficients one by one, but suppress the detailed calculations. The latter
can be carried out by hand, though the procedure itself was refined by
computer (Maple and Mathematica versions are available from the authors).

The system defined by \eqref{generic} does not in general define a Lie algebra
as $d^2$ (it is clearer to type this as $dd$) may not vanish. A valid solution
must therefore satisfy both $dd\a i=0$ and the SKT condition $\sq\Om=0$. On a
computer, the operator $\sq=-\ft12d\cir(-iJ)\cir d$ can be executed by
treating $-iJ$ as the substitution $\a\j\mapsto-\a\j$. If $dd=0$ then
$\sq=\pd\opd$, but in general $\sq\Om$ may have a non-zero (3,1) component. To
avoid this embarassment, we shall restrict attention to
\begin{equation}\label{ijOm}\aa i\j\sq\Om=\a i\we \a\j\we\sq\Om\end{equation}
which we identify with its coefficient relative to the standard 6-form
$\a{123\1\2\3}$. (Actually we shall only ever take $i=j$ in \eqref{ijOm}.)

A first calculation reveals that $\ao11\sq\Om=\xo33|c_1|^2$. The positive
definiteness of $\Om$ implies that $c_1=0$. Independently, the $\a{23\2}$
component of $dd\a3$ equals $|c_4|^2$, so $c_4=0$. To simplify matters, we now
consider two cases: $a_5\ne0$ and $a_5=0$ (the choice of subscript $5$ is
partly a matter of taste).\vs

\n\emph{Case 1.} Suppose $a_5\ne0$. Subtracting a multiple of $\a2$ from
$\a3$ we may suppose that $b_5=0$. Using $dd\a3$ (we mean of course `the
vanishing of $dd\a3$') gives $c_3=0$. Using $dd\a2$ now gives $a_4=0$ and
$b_1=b_2=0$. To sum up, $a_4,b_1,b_2,b_5,c_1,c_3,c_4$ are all zero. If $p$
denotes the $\a{13\1}$ component of $dd\a3$ then $p=a_2c_2+b_4\ol a_5$ and
\[ 0=\ao22\sq\Om=(|a_2|^2+|a_5|^2)\xo22-(\Re\,p)\xo33.\]
This contradicts $a_5\ne0$. $\frown$\vs

\n\emph{Case 2.} Suppose $a_5=0$. We divide into two subcases: $c_3\ne0$ and
$c_3=0$.\vs

\n\emph{Case 2.1.} Suppose $c_3\ne0$. Inspecting only $dd\a3$ already gives
$a_1=a_2=0$ and $b_5=0$. Using $dd\a2$ then gives $a_4=0$. Returning to $dd\w3$
gives $b_2=a_3=0$. To sum up, $a_*,b_2,b_5,c_1,c_4$ are all zero so we recover
\eqref{SKT1}. The proof of Lemma~\ref{vanish} amounts to nothing more than an
application of the equation $\ao33\sq\Om=0$, which yields \eqref{SKT2}. In this
way we arrive at the conclusion of Theorem~\ref{SKT} that represents the
generic solution to the SKT hypothesis. $\smile$\vs

\n\emph{Case 2.2.} Suppose $c_3=0$. Using $dd\a2$ shows that $a_4=0$ (even if
$a_1=a_2=0$). We now divide into two subsubcases: $c_2\ne0$ and $c_2=0$.\vs

\n\emph{Case 2.2.1.} Suppose $c_2\ne0$. Using $dd\a2$ and $dd\a3$ respectively
gives $a_2=0$ and $b_5=0$. Then $\ao22\sq\Om =\xo33|b_2|^2$, whence
$b_2=0$. Using $dd\a3$ gives $a_1=0$. It follows that
\[0=\ao33\sq\Om=\xo33(|b_1|^2+|b_4|^2+|c_2|^2)\] and the solution reduces to
\begin{equation}\label{subsumed}\left\{\ba{rcl}d\a1&=&0\\d\a2&=&a_3\wo11\\
d\a3&=&b_3\ao11.\ea \right.\end{equation} By subtracting a multiple of $\a3$
from $\a2$ (or swapping the two if $b_3=0$), this solution can be subsumed
into that of Theorem~\ref{SKT}. $\smile$\vs

\n\emph{Case 2.2.2.} Suppose $c_2=0$ so that $c_*=0$. The vanishing of $dd\a3$
implies that $b_5=0$ and consequently that either $a_1=a_2=0$ or $b_4=0$. In
the former case, the vanishing of $\ao22\sq\Om$ and $\ao33\sq\Om$ gives
$b_2=0$ and $b_1=b_4=0$ respectively, and we obtain \eqref{subsumed}. The
final situation to deal with is therefore that $a_4,a_5,b_4,b_5,c_*$ all
vanish. This implies that \begin{equation}\label{complexline}\ba{rcl} 0\ =\
\ao22\sq\Om &=& \xo22|a_2|^2+\xo23a_2\ol b_2+\xo32\ol a_2b_2+\xo33|b_2|^2\y
&=&\Om(a_2V_2+b_2V_3,\,\ol a_2\ol V_2+\ol b_2\ol V_3),\ea\end{equation} in the
notation of \eqref{restriction}. Unless $a_2=b_2=0$, the restriction of $\Om$
to the complex line spanned by $a_2V_2+b_2V_3$ is zero, which is
impossible. Exactly the same argument applied to $\ao33\sq\Om$ gives
$a_1=b_1=0$ and we are left with only $a_3,b_3$ non-zero, whence
\eqref{subsumed}. $\smile$

\begin{oss} An invariant complex structure $J$ always induces a complex Lie 
algebra structure on the $i$-eigenspace $\g^{1,0}$ of $\g_c$. In the
6-dimensional nilpotent case, $\g^{1,0}$ is either abelian or isomorphic to the
complex Heisenberg algebra. In the former case, $J$ is itself called
\emph{abelian}, and this is equivalent to asserting that $d$ maps the
subspace $\La^{1,0}$ of $\g^*_c$ into $\La^{1,1}$. The complex structure given
by Theorem~\ref{SKT} is therefore abelian if and only if $E=0$. So SKT does not
imply that the complex structure is abelian. This is in contrast with the
result that if a $2$-step nilpotent Lie group admits an invariant HKT structure
then the hypercomplex structure must be abelian \cite{DoFHTS}.\end{oss}

\section{Real Lie algebras}

Theorem~\ref{SKT} can be used to classify explicitly the real $6$-dimensional
nilpotent Lie algebras admitting an SKT structure. Before explaining this, we
shall introduce a formalism that will eventually help to understand and
manipulate \eqref{SKT1}. Namely, we shall identify forms of type (1,1) with
$2\times2$ matrices by setting
\begin{equation}\label{Yal}\Y_\al=A\oa12+B\oa22+C\ao11+D\ao12\end{equation}
where
\begin{equation}\label{Y}\Y=\Big(\!\ba{cc}A&B\\C&D\ea\!\Big)\end{equation} so
that \begin{equation}\label{da3Y} d\a3= \Y_\al+E\aa12.\end{equation} The exact
positioning of the coefficients may seem strange, but follows a logic that is
revealed in \eqref{minors} below. The subscript $_\al$ indicates the basis
relative to which the construction is made.

Operations on 2-forms now translate into matrix operations in a natural way.
For example, \[\X_\al\we\Y_\al=\tr(\X\Y^\#)\a{1\12\2}\] where \[\Y^\#=\hbox{
adj}\,\Y=\Big(\!\ba{cc}D&-B\\-C&A\ea\!\Big)=(\det\Y)\Y^{-1}\] is the transpose
of the matrix of cofactors. Moreover, \[\ol{\Y_\al}=\ol\Y^\#\sb\al\] and so
\[\Y_\al\we\ol{\Y_\al}=\tr(\Y\ol\Y)\a{1\12\2}.\]

\begin{exa} The complex coefficients $A,B,C,D$ are meant to be thought of as
those in \eqref{SKT1}. In an illustration, we compute
\begin{equation}\label{triple}\ba{rcl} d\a3\we d\a\3 &=&
(\tr(\Y\ol\Y)-|E|^2)\a{1\12\2}\y Jd\a3\we d\a\3 &=&
(\tr(\Y\ol\Y)+|E|^2)\a{1\12\2}\y d\a3\we d\a3 &=&
\tr(\Y\Y^\#)\a{1\12\2}=2(\det\Y)\a{1\12\2}.\ea\end{equation} The middle
equation is a more succinct version of the SKT formula \eqref{ABC}.\end{exa}

The following result relies on the classification of \cite{SalCSN}, whose
notation we freely adopt.
 
\begin{teo}\label{real} A $6$-dimensional nilmanifold $M=\Ga\bs G$
admits an invariant SKT structure if and only if the Lie algebra $\g$
is isomorphic to one of \[\begin{split} &(0,0,0,0,13+42,14+23)\\
&(0,0,0,0,12,14+23)\\ &(0,0,0,0,12,34)\\ &(0,0,0,0,0,12).\end{split}\]
In particular, $\g$ is 2-step and has first Betti number
$b_1(\g)=b_1(M)$ at least 4.\end{teo}

\begin{D} First note that the possibility that $\g$ is abelian is precluded by
Definition~\ref{strong}. The vanishing the the real and imaginary components of
$d\a1,d\a2$ in \eqref{SKT1} implies immediately that $b_1(\g)\ge4$. The fact
that $d\a3\in\alt^2(\ker d)$ means (in the notation of \cite{SalCSN}) that
$(\g^2)^\mathrm{o}=V_2$ equals $\g$, which is therefore 2-step.

Using the methods of \cite{SalCSN}, any 2-step nilpotent Lie algebra with
$b_1\ge4$ is isomorphic to one of \begin{equation}\label{list}\ba{llll}
\mathrm{(i)}&(0,0,0,0,12,13)&\mathrm{(ii)}&(0,0,0,0,13+42,14+23),\y
\mathrm{(iii)}&(0,0,0,0,12,14+23)\qq&\mathrm{(iv)}\ &(0,0,0,0,12,34)\y
\mathrm{(v)}&(0,0,0,0,0,12)&\mathrm{(vi)}&(0,0,0,0,0,12+34)\ea\end{equation}
For example, in case (i) there is a real basis $(e^i)$ of 1-forms for which
$de^i=0$ for $1\le i\le4$, $de^5=e^{12}$ and $de^6=e^{13}$. We need to
eliminate (i) and (vi), and prove existence in the other cases.

Given \eqref{SKT1}, write $d\a3=\si^1+i\si^2$ in real and imaginary
components, and consider the real $2\times2$ matrix $\B=(b^{ij})$ associated
to the bilinear form \[\si^i\we\si^j=b^{ij}\a{1\12\2}.\] Under the SKT
assumption, equations \eqref{triple} give
\begin{equation}\label{B}-\B=\left(\!\ba{cc}|E|^2+\Re U&\Im U\\\Im U
&|E|^2-\Re U\ea\!\right)\end{equation} where $U=-\det\Y=BC-AD$. Using
$de^5,de^6$ in place of $\si^1,\si^2$ is merely a change of real basis and
must yield a matrix congruent to $\B$. It follows that, in the above
examples,\vs

$\B$ is the zero matrix for (i) and (v)

$\B$ has rank 1 for (iii) and (vi)

$\det\B\ne0$ for (ii) and (iv).\vs

\n In case (vi), we may rescale $\a3$ so that $d\a3$ is real. This implies
that $d\a3$ is a (1,1) form and $E=0$. Since $\B$ has rank 1, the matrix
\eqref{B} has zero determinant so $U=0$. This means that $\B=0$, which is a
contradiction. In case (i) we already know that $\B=0$ so that $E=0=U$. Thus,
$d\a3=\Y_\al$ is a (1,1) form with $\rank\,\Y\le1$. But the image $\langle
e^{12},e^{13}\rangle$ of $d$ in $\alt^2g^*$ is divisible by the real 1-form
$e^1$ and must therefore be generated by $e^1\we Je^1$. But this contradicts
the fact that $d(\g^*)$ is actually 2-dimensional.

The remaining cases are distinguished by the rank and signature of $\B$, and it
is easy to check that the coefficients in \eqref{SKT2} can be chosen to realize
the four possibilities.\end{D}

\begin{exa} The irrelevance of the choice of Hermitian metric is special to
the nilpotent situation. The third Lie algebra listed in Theorem~\ref{real}
corresponds to the product $H_3\times H_3$ where $H_3$ is the real Heisenberg
group.  A simple example in which the SKT condition is metric dependent is
provided by $H_3\times S^3$, where $S^3$ is identified with $SU(2)$. We may
choose a real basis of 1-forms such that \[de^1=0,\q de^2=0,\q de^3=e^{12},\q
de^4=e^{56},\q de^5=e^{64},\q de^6=e^{45}.\] Setting
$\w1=e^1+ie^2,\,\w2=e^3+ie^4,\,\w1=e^5+ie^6$ gives \[\left\{\ba{l} d\w1=0\y
d\w2=\ft12(i\,\wo11-\wo33)\y d\w3=\ft12(\ww23+\wo32)\ea\right.\] It follows
that $d(\wo33)=0$ and $\pd\opd\wo22=d\w2\we\ol{d\w2}=0$, so \eqref{Om}
satisfes $\sq\Om_0=0$.  On the other hand, $\sq\wo13=-\ft14i\w{12\2\3}$ and
the general 2-form \eqref{Omega} determines an SKT metric if and only if
$\xo13=-\ol{\xo31}=0$.\end{exa}

\section{Invariant forms on Iwasawa}

We first summarize the relevant facts concerning left-invariant complex
structures on the Iwasawa manifold. We have attempted to give a self-contained
account, though important background for \S\S4,5 can be found in
\cite{KeSCSI,SalCSN}. The reader is implicitly referred to these papers for
further explanation of a number of points.

Let \[G=\left\{\left(\!\ba{ccc}1&z^1&z^3\\0&1&z^2\\0&0&1\ea\!\right):z^i\in\C
\right\}\] denote the complex Heisenberg group and $\Ga$ the discrete subgroup
for which $z^i$ are Gaussian integers. We define $M$ to be the set
$\Ga\backslash G=\{\Ga g:g\in G\}$ of right cosets. It is a homogeneous space
relative to the action of $G$ by right translation that persists on the
quotient, though we shall be interested in the projections of tensors that are
invariant by left translation on $G$.

The complex 1-forms $\w1=\>dz^1,\,\w2=dz^3,\,\w3=-dz^3+z^1dz^2$ satisfy
$d\w3=\ww12$ and span the the $(1,0)$ space $\La_0$ of the bi-invariant
complex structure $J_0$ on $M$. It is known that, in addition to $J_0$, any
left-invariant complex structure on $G$ leaves invariant the real
4-dimensional subspace \begin{equation}\label{DD}\DD=
\langle\Re\w1,\Re\w1,\Im\w2,\Im\w2\rangle\end{equation} that arises from a
principal $T^2$-fibration $M\to T^4$. It therefore makes sense to consider the
space $\La=\langle \a1,\a2,\a3 \rangle$ generated by the modified 1-forms
\begin{equation}\label{aaa}\left\{\ba{l}\a1=\w1+a\w\1+b\w\2\\
\a2=\w2+c\w\1+d\w\2\\\a3=\w3+x\w\1+y\w\2+u\w\3\ea \right.\end{equation} with
$a,b,c,d,x,y,u\in\C$. If $\La$ is maximally complex it defines an invariant
almost complex structure on $M$ that we denote by $J_{\X,x,y}$.

The effect of $J_{\X,x,y}$ on a real basis $(e^i)$ can be deduced by setting
\[\w1=e^1+ie^2,\q\w2=e^3+ie^4,\q\w4=e^5+ie^6,\] though little is to be gained
from this. The integrability condition for $J_{\X,x,y}$ is readily expressed
in terms of \eqref{aaa} as $d\a3\we\a{123}=0$ (equivalently $d\a3\we\aa12=0$).
This reduces to the equation \begin{equation}\label{u}
u=bc-ad=-\det\X\end{equation} where \[\X=\Big(\!\ba{cc}a&b\\c&d\ea\!\Big).\]
Using \eqref{Yal}, we may consider the simple 2-form
\begin{equation}\label{minors}\ba{rcl}\aa12
&=&\ww12+a\ow12+b\ow22+c\wo11+d\wo12-u\ww\1\2\y
&=&\ww12+c\wo11+d\wo12-a\wo21-b\wo22-u\ww\1\2\y &=&\X_\om+\ww12-u\ww\1\2.\ea
\end{equation} We write the characteristic polynomial of $\XoX$ as
$c(x)=x^2-\ga x+\de$, so that
\begin{equation}\label{ga}\ga=\tr(\XoX),\qq\de=\det(\XoX)=|u|^2\end{equation}
(notation of \cite{KeSCSI}).

The formulae
\begin{equation}\label{vols}\ba{l}\a{1\12\2}=-\aa12\we\aa\1\2=(1-\ga+\de)
\w{1\12\2}=c(1) \w{1\12\2}\y
\a{1\12\23\3}=\a{1\12\2}\we\ao33=c(1)(1-\de)\w{1\12\23\3}\ea\end{equation}
express volume changes associated to a switch of basis from $\om$ to $\al$. As
a consequence, $\La\cap\ol\La=\{0\}$ if and only if
\begin{equation}\label{c(1)}\de\ne1\q\hbox{ and }\q c(1)\ne0\end{equation} and
these are the conditions that ensure that $J_{\X,x,y}$ is well defined. From
now on, we assume that \eqref{u} and \eqref{c(1)} hold. For simplicity, we
also suppose $x=y=0$, and denote $J_{\X,0,0}$ simply by $J_\X$. It will become
obvious that reducing to this case causes no loss of generality, the key point
being that $d\a3$ (and so \eqref{sharp} below) does not involve $x,y$. An
underlying reason for the irrelevance of $x,y$ is provided in \S5.

One may invert \eqref{minors} so as to express $\ww12$ in terms of the
$\al$'s. Up to the overall factor $c(1)$, the corresponding formula is given 
by reversing the signs of $a,b,c,d$:

\begin{lem}\label{reverse} $c(1)\ww12=-\X_\al+\aa12-u\aa\1\2$.\end{lem}

\begin{D} Consider the bases $(\w1,\w2,\w\1,\w\2),(\a1,\a2,\a\1,\a\2)$. The
second is related to the first by the block matrix \[\M=\left(\!\ba{c|c}
\I&\X\\\hline\vph\ol\X&\I\ea\!\right)\] in a row-by-row fashion. Set
$\Z=(\I-\XoX)^{-1}$ so that $c(1)=\det(\Z^{-1})$. By first observing that
\[\X\ol\Z=\Z(\Z^{-1}\X) \ol\Z=\Z(\I-\XoX)\X\ol\Z
=\Z\X(\I-\ol\X\X)\ol\Z=\Z\X,\] it is easy to verify that
\[\M^{-1}=\left(\ba{c|c}\Z&-\Z\X\\\hline\vph -\ol\X\Z& \ol\Z\ea\right).\]

The coefficients $1,c,d,-a,-b,-u$ featuring in \eqref{minors} are the
$2\times2$ minors of the $2\times4$ matrix $\M=(\I\,|\,\X)$, corresponding to
the Mathematica command $\Minors[\M,2]$. The coefficients of $\ww12$ are
therefore given by \[\Minors[(\Z\,|-\!\!\Z\X),\,2]=(\det \Z)
\Minors[(\I\,|-\!\X),\,2],\] using a well-known property of minors.\end{D}

\begin{pro}\label{lim} The complex structure $J_X$ is SKT if and only if the
eigenvalues of $\XoX$ satisfy the equation
\begin{equation}\label{curve}(1+|z|^2)\,|1\!+\!z|^2=8|z|^2\end{equation}
illustrated below.\end{pro}

\n It follows easily that all points of the curve are realized except for $z=1$
that is excluded by \eqref{c(1)} (see Theorem~\ref{main} below).

\vspace{0pt}\hspace{70pt}\scalebox{0.6}{\includegraphics{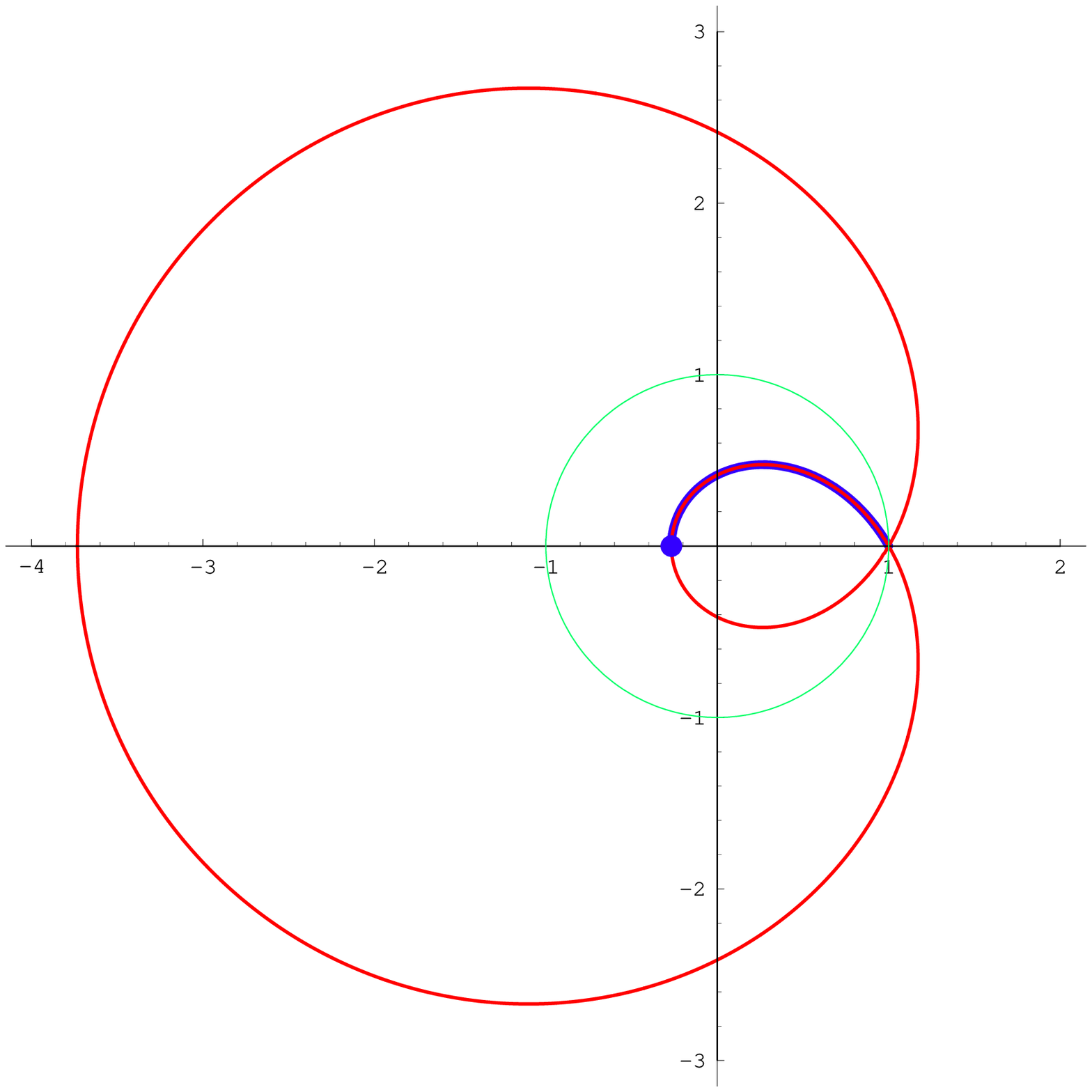}}
\vspace{20pt}

\begin{D} We first prove that $J_\X$ is SKT if and only if
\begin{equation}\label{de}1-6\de+ \de^2+\ga+\ga\de=0.\end{equation} Using
Lemma~\ref{reverse}, \begin{equation}\label{sharp}\ba{rcl}c(1)d\a3 &=&c(1)
\ww12+uc(1) \ww\1\2\y &=&-\X_\al+\aa12-u\aa\1\2+u(-\ol\X^\#\sb\al+ \aa\1\2-\ol
u\aa12)\y &=&-(\X_\al+u\ol\X^\#\sb\al)+(1-\de)\aa12\y &=&
(-\X+\de\ol\X^{-1})_\al+(1-\de)\aa12.\ea\end{equation} Comparing this with
\eqref{da3Y} yields \begin{equation}\label{YE} c(1)\Y=-\X+\de\ol\X^{-1},\qq
c(1)E=1-\de.\end{equation} Thus
\[c(1)^2\tr(\Y\ol\Y)=\tr(\XoX-2\de\I+\de^2(\XoX)^{-1})=\ga-4\de+\de^2
(\ga/\de)\] and \eqref{de} follows from Lemma~\ref{vanish} and \eqref{triple}.

Now let $\la,\mu$ denote the eigenvalues of $\XoX$, so that $\ga=\la+\mu$ and
$\de=\la\mu$. Thus, $\la\mu\ge0$ so $\la$ and $\mu$ cannot be real with
opposite signs. It is an elementary but non-trivial fact that if $\la,\mu$ are
real and non-positive then they are equal \cite[Problem~6,\,\S4.6]{HoJMaA}.
Under the SKT assumption, \[\ga=\frac{-1+6\de-\de^2}{1+\de}\] is non-negative
if and only if $3-2\sqrt2\le\de\le3+2\sqrt2$. The eigenvalues $\la,\mu$ can
only be real if
\[0\le(\la-\mu)^2=\ga^2-4\de=\left(\frac{1-\de}{1+\de}\right)^2(\de^2-14\de+1)
\] which implies that $\de\le7-\sqrt{48}$ or $\de\ge7+\sqrt{48}$. The various
inequalities are incompatible, and there are no solutions with $\la,\mu>0$ and
it follows that $\la=\ol\mu$. Thus, \[1-\ga+\de=1+\la+\ol\la+|\la|^2=|1+\la|^2
>0\] and \eqref{de} translates into \eqref{curve}.\end{D}

\begin{oss}\label{inversion} It is easy to check that equation \eqref{curve} 
is unchanged by substituting $1/z$ in place of $z$. We shall see below that
this corresponds to reversing the sign of $J$. Another curve with a similar
shape invariant by $z\mapsto1/z$ and $z\mapsto\ol z$ is given in polar
coordinates by $r=e^{\sin\th}$. (This can be generalized by replacing $\sin\th$
by an odd Fourier series.) If we shift the origin to the point $z=1$ of
self-intersection, \eqref{curve} becomes \[r=-3\cos\th\pm\sqrt{2+\cos^2\th}.\]
A true lima\c con has the somewhat simpler equation $r=2\cos\th+1$, but the
corresponding translate is not invariant under inversion.\end{oss}

\section{Moduli space interpretation}

Let $\cC$ denote the set of all invariant complex structures on the Iwasawa
manifold $M$. This is a subset of the set of all almost complex structures on
$\g$ (equivalently, maximally complex subspaces of $\g_c$), that can in turn be
identified with the homogeneous space $GL(6,\R)/GL(3,\C)$. It is known that
$\cC$ has four connected components, and these can be described as follows.

Changing the sign of $J$ correponds to an overall reversal of orientation, and
corresponds to the
transformation\begin{equation}\label{inv}(\a1,\a2,\a3)\mapsto(\a\1,\a\2,\a\3).
\end{equation} This identifies the components of $\cC$ in pairs. Let $\hJ$
denote the restriction of $J$ to the real subspace \eqref{DD} underlying
$\langle\w1,\w2 \rangle$. The remaining two components of $\cC$ are
distinguished by the orientation of $\hJ$ or equivalently, by \eqref{vols},
the sign of $c(1)$.  For example
\begin{equation}\label{cC+}\cC_+=\{J_{\X,\,x,y}: x,y\in\C,\ |u|<1,\
c(1)>0\}\end{equation} is the connected component that contains $J_0$.

To fully understand the SKT constraint on complex structures on $M$, we need to
describe various group actions on $\cC$ and their effect on the matrix $\X$. 
\vs\vs

\n\textbf{\n(i) Involution.} Referring to \eqref{aaa}, can say that the almost
complex structure $\hJ$ is represented by the matrix $(\I\,|\,\X)$. This is
replaced by \begin{equation}\label{trans}(\ol\X\,|\,\I)
\equiv(\I\,|\,\ol\X^{-1})\end{equation} under \eqref{inv} when $J,\hJ$ are
replaced by $-J,-\hJ$. The equivalence relation $\equiv$ indicates
premultiplication on both halves by an invertible matrix, and reflects a
re-adjustment \[(\a\1,\a\2)\mapsto(\ti\al^1,\ti\al^2) \] into the row echelon
form \eqref{aaa} in which $\ti\al^i$ has leading term $\w i$. Notice that
$\det(\ol\X^{-1})=-1/\ol u$, consistent with setting $\ti\al^3=\a\3/\ol
u=\w3+\w\3/\ol u$. It now follows that \[-(J_\X)= J_{\ol\X^{-1}}.\]

If the eigenvalues of $\X$ are $\la,\ol\la$ then those of $\ol\X^{-1}$ are
$1/\ol\la,\,1/\la$, and this justifies Remark~\ref{inversion} and helps to
explain why the solution curve has a lima\c con shape. The inner part $|z|<1$
corresponds to solutions in the same component as $J_0$ and the outer part
$|z|>1$ to the component of $-J_0$. One is an inversion of the other. The
thicker part with $|z|<1$ and $\Im z\ge0$ fits (after rotation by
$45^\mathrm{o}$) into the semicircular region of the diagram in \cite{KeSCSI},
and represents the SKT solutions most faithfully. (The rest of the complex
plane above does not fit into the diagram that is a schematic pasting of the
real and complex plane best contemplated in 3 dimensions.)\vs

The integrability condition for a left-invariant almost complex structure on
an arbitrary Lie algebra $\g$ is invariant under the action of the
automorphism group \[\cG=\{ f\in\End\g:[f(v),f(w)]=f[v,w].\}\] In the Iwasawa
case, $\cG$ can be identified with the semidirect product
$GL(2,\C)\ltimes\C^2$ consisting of complex $3\times3$ matrices
\begin{equation}\label{semi}\left(\ba{c|c}\vph\q\P\q&\mathbf{q}\yy \hline\vph
0&\det\P\ea\!\!\right),\qq\P\in GL(2,\C),\q\mathbf{q}\in\C^2\end{equation} (see
\cite{SalCSN}). The action on $\cC$ of $\P$ and $\mathbf q$ can be considered
separately.\vs\vs

\n\textbf{(ii) Right translation.} The normal subgroup $\C^2$ of $\cG$ can be
identified with the group $\Ad G$ of inner automorphisms of $\g$. Since
$J\in\cC$ is by definition left-invariant, $\Ad(g)J$ equals the right translate
of $J$ by $g^{-1}\in G$. Whilst $J_0$ is fixed by this action ($G$ being a
\emph{complex} Lie group), all the other orbits have positive
dimension. Moreover, if $u\ne0$ then $J_{\X,x,y}$ lies in the same orbit as
$J_\X$ for any $x,y\in \C$. This is explained in \cite{KeSCSI}, but to avoid
duplication we next insert an infinitesimal version of this fact.

\begin{oss} The action of the diffeomorphism group on a complex structure 
$J$ is detected by the image of \[\opd\colon\Ga(T\otimes\La^{0,0})\to
\Ga(T\otimes\La^{0,1})\] where $T$ denotes the holomorphic tangent bundle of
$(M,J)$. Restricting to invariant tensors, a $p$-form $f$ with values in $T$
can be regarded as a linear mapping $\La^{1,0}\to\La^{p,0}$ and $\opd f$ is
calculated by means of the formula \[(\opd f)(\a i)=\opd(f(\a i))-f(\opd\a i)=
-f(\opd\a i).\] Since $\opd\a3=\Y_\al$, the dimension of $\Im\opd$ coincides
with the rank of $\Y$ which equals 2 if and only if $u\ne0$ (see \eqref{da3Y}
and \eqref{YE}). This discussion is relevant to the jumping of Hodge numbers
discussed in \cite{NakCPM}.\end{oss}

In the SKT context, it is no restriction to impose the condition $\det\X\ne0$
since $z=0$ is not a solution of \eqref{curve}. In the light of the above
remarks, one might call the structures $J_{\X,x,y}$ for which $\det\X\ne0$ the
\emph{stable} points of $\cC$. Since $G$ acts smoothly on $M$ by right
translation, $J_\X$ and $J_{\X,x,y}$ then determine the same point in the
moduli space of complex structures modulo diffeomorphism, and for calculations
we may assume that $x=y=0$ so that $J=J_\X$ is completely determined by $\hJ$
and the matrix $\X$.\vs

\n\textbf{(iii) Outer automorphisms.} The quotient $\cG/\Ad(G)$ can be
identified with $GL(2,\C)$, an element of which acts by a change of basis
\[\left\{\ba{l}\w1\mapsto p^1_1\w1+p^1_2\w2\y\w2\mapsto p^2_1\w1+p^2_2
\w2\ea\right.\] By analogy to \eqref{trans}, the matrix $(\I\,|\,\X)$
representing $\hJ$ is transformed into \[(\P\,|\,\X\ol
\P)\equiv(\I\,|\,\P^{-1}\X\ol \P)\] and therefore a left action on $\cC$ is
defined by \[\P^{-1}\cdot J=J_{\P^{-1}\X\ol \P}.\] The presence of $\det\P$ in
\eqref{semi} ensures that the extension from $\DD$ to $\g$ is well defined.

The remaining action by $GL(2,\C)$ is less geometrical and the resulting
quotient $\cC/\cG$ has singularities, an example of which is given after
Example~\ref{sing}. Now
\begin{equation}\label{consim}\X\mapsto\P^{-1}\X\ol\P,\qq\P\in
GL(n,\C)\end{equation} is an action that gives rise to the theory of
\emph{consimilarity} for $n\times n$ matrices $\X$. If $[\![\X]\!]$ denotes
the consimilarity class of $\X$ (i.e.\ an orbit for the above action) and
$[\Y]$ the similarity class of $\Y$, there is a well defined mapping
\begin{equation}\label{orbits}\phi\colon[\![\X]\!]\mapsto[\XoX].\end{equation}
This is not a bijection as $\XoX$ can be zero without $\X$ being
zero. However, the general theory of consimilarity developed in
\cite[\S4.6]{HoJMaA} implies that $\phi$ restricts to a bijection between
classes subject to an equal rank condition. In the simple case of $n=2$, we
may condense the discussion of this section into

\begin{teo}\label{main} An invariant complex structure $J$ on the Iwasawa 
manifold is SKT if and only if it equals $J_{\X,x,y}$ where
\[\X=\P^{-1}\Big(\!\ba{cc}0&z\\ 1&0\ea \!\Big)\ol\P\] with $\P\in GL(2,\C)$,
$z$ a solution of \eqref{curve} different from 1 and
$x,y\in\C$. Moreover, such structures lie in the connected components
of $J_0$ and $-J_0$ in $\cC$.\end{teo}

\begin{D} If $\X$ has the form given then $\XoX$ is similar to a diagonal
matrix with entries $z,\ol z$. Since $z\ne1$, $J_{\X,x,y}$ is well defined.
Bearing in mind (from (ii)) that $x,y$ are irrelevant, Proposition~\ref{lim}
implies that $J_{\X,x,y}$ is SKT.

Conversely, suppose that $J\in\cC$ is SKT. The complex structures $J_{\X,x,y}$
represented by \eqref{aaa} constitute a dense set of $\cC$, and the only
missing points are those arising when one or more of the coefficents become
infinite. The SKT condition involves only $d\a3$, so the only case potentially
not covered by \eqref{aaa} is that in which $\a3$ belongs to the span of $\w\3$
and $\DD$. But then $Jd\a3\we d\a\3\ne0$ and $J$ cannot be SKT by
Lemma~\ref{vanish}. The statement about connected components now follows from
the fact that the solutions in Proposition~\ref{lim} all satisfy $c(1)>0$ (see
\eqref{cC+}).

We may therefore suppose that $J=J_{\X,0,0}=J_\X$, and that the eigenvalues
$z,\ol z$ of $\XoX$ satisfy \eqref{curve}. Assume firstly that $z$ is not
real. Choose $\P\in GL(2,\C)$ such that $\Y=\P\X\ol\P^{-1}$ satisfies
\[\Y\ol\Y=\P(\XoX)\P^{-1}=\Big(\!\ba{cc}z&0\\0&\ol z\ea\!\Big).\] In the
notation \eqref{Y}, $A\ol B+B\ol D=0= A\ol C+C\ol D$, which implies \[A(B\ol
C-\ol BC)=0 =D(B\ol C-\ol BC),\] and $A=D=0$. We can now premultiply $\P$ by a
diagonal matrix so as to convert $B$ to $z$ and $C$ to 1, as required.

Now assume that $z\in\R$. The equation $\Y\ol\Y=\Big(\!\ba{cc}z&0\\
1&z\ea\!\Big)$ would imply that $B=0$ and $z>0$, which contradicts
\eqref{curve}. Thus, $\XoX$ is again diagonalizable and we need to solve
\begin{equation}\label{YY}\Y\ol\Y=z\,\I,\end{equation} that implies that
$|A|=|D|$. But we can find $\Q\in GL(2,\C)$ such that $\ti\Y=\Q\Y\ol\Q^{-1}$
has its last entry $\ti D$ zero. Since $\ti\Y$ satisfies \eqref{YY} in place
of $\Y$, we also have $\ti A=0$, and we can modify $\Q$ so that $\ti C=1$ and
$\ti\Y=\Big(\!\ba{cc}0&z\\ 1&0\ea\!\Big)$.\end{D}

\begin{exa}\label{sing} The point $z=1$ on the curve is not admissible,
because in this case \[\a1=\w1+\w\2= \ol{\a2},\] and any corresponding tensor
$J$ is degenerate. The solutions \[z= -2\pm\sqrt3\in\R\] are admissible, but
since $-2-\sqrt3=1/(-2+\sqrt3)$ the $J$'s coincide up to an overall complex
conjugation, in accordance with the discussion in (i). These are the simplest
solutions to the SKT equation. The purely imaginary solutions are $z=\pm
i\sqrt{3\pm2\sqrt2}$, and we presume that there are no rational solutions to
\eqref{curve}.\end{exa}

Let $H_z$ denote the stabilizer of $\Big(\!\ba{cc}0&z\\1&0\ea\!\Big)$ in
$GL(2,\C)$ for the action \eqref{consim}. It is easy to verify that if
$z\in\R$ and $z\ne0$ then
\begin{equation}\label{Hz}H_z=\left\{\Big(\!\ba{cc}a&-\ol cx\\c&\ol
a\ea\!\Big): |a|^2+|b|^2\ne0\right\}\cong GL(1,\mathbb{H})\end{equation} If
$z\in\C\setminus\R$ then $H_z$ is isomorphic to the common subgroup $\C^*$
obtained by setting $c=0$ in \eqref{Hz}. It follows that the SKT structures
$z=-2\pm\sqrt3$ (one of which is a blob in the Figure) represent singular
points in $\cC/\cG$.

\section{Balanced structures and reduced holonomy}

Let $M$ be a nilmanifold of dimension $2n$ with an invariant complex
structure. The basis arising from Theorem~\ref{teosal} furnishes us with a
closed form \begin{equation}\label{eta}\eta=\a{12\cdots n}\end{equation} of
type $(n,0)$. Being invariant, $\eta$ also has constant norm, so it is natural
to ask whether it is parallel with respect to a suitable connection.

Throughout this final section, $\na$ denotes any Hermitian connection,
that is one satisfying $\na g=0$ and $\na J=0$. We use $D$ to denote
the Levi-Civita connection, so that $Dg=0$ but in general
$DJ\ne0$. The holonomy of $\na$ is contained in $SU(n)$ if and only if
$\na\eta=0$.

\begin{pro}\label{holo} Suppose that $(J,g)$ is an invariant Hermitian 
structure on a nilmanifold $M$ of dimension $2n$. Then the holonomy of $\na$
reduces to $SU(n)$ iff $g$ is balanced.\end{pro}

\begin{D} Since $\na$ preserves the complex structure, we may write \[\na\eta=
i\beta\otimes\eta.\] The connection 1-form $i\beta$ belongs to the Lie algebra
$\mathfrak{u}(1)$ and $\beta$ is \emph{real} (this can be checked directly by
applying $\na$ to $\|\eta\|^2=\langle\eta,\ol\eta\rangle$).

The tensor $DJ$ can be identified with both $D\Om$ and $d\Om$, and gets
converted into the torsion of $\na$ in passing from $D$ to $\na$. It follows
that \[\na\eta=D\eta+d\Om\cdot\eta\] where $\>\cdot\>$ stands for a suitable
linear mapping. If we skew-symmetrize both sides, the Levi-Civita connection
$D$ gets converted to $d$, and extracting $(n,1)$ components gives
\begin{equation}\label{beta}\beta^{0,1}\we\eta=d\eta+(d\Om
\cdot\eta)^{n,1}=(\opd\Om\cdot\eta)^{n,1}.\end{equation} The last expression
determines a $U(n)$ equivariant mapping
\begin{equation}\label{contraction}\La^{1,2}\otimes\La^{n,0}\to
\La^{n,1}\cong\La^{n-1,0}\end{equation} and the (2,1)-component $\pd\Om$ is
omitted because it can only define the zero map $\La^{n,0}\to \La^{n,1}$. But
\eqref{contraction} is equivalent to the $SU(n)$-equivariant contraction
$\La^{1,2}\to\La^{0,1}\cong\La^{n-1,0}$ that extracts the `trace' of $\opd
\Om$. It follows that $\eta$ is parallel if and only if $\th=0$.\end{D}

\n Combined with Proposition~\ref{incomp}, this is consistent with a result of
\cite{PapKHG}, namely that a compact $2n$-dimensional manifold endowed with a
conformally balanced Hermitian structure for which $\sq\Om=0$ and $\na$ has
holonomy in $SU(n)$ is in fact K\"ahler and therefore Calabi-Yau.\vs

Returning to six dimensions, it is an easy matter to list Lie algebras
admitting a balanced KT structure for which the holonomy of $\na$ therefore
reduces. We restrict attention to those considered in \S3. Case (v) is realized
by taking $C=1$ and $A,B,D,E$ zero which is incompatible with \eqref{bal}. But
this is the only one excluded:

\begin{cor} Each of the Lie algebras \[\begin{split}&(0,0,0,0,13+42,14+23)\\
&(0,0,0,0,12,14+23)\\ &(0,0,0,0,12,34)\\ &(0,0,0,0,0,12,13)\\
&(0,0,0,0,0,12+34)\end{split}\] admits a Hermitian structure for which
$\na$ has $SU(3)$ holonomy.\end{cor}

\begin{D} It suffices to consider the standard metric corresponding to 
\eqref{Om}. The `balanced' condition $B=C$ is then independent of $A,D,E$, and
this gives us the flexibility to realize the various cases as follows:

\n If $B=C=0$ then

$A=D=0$ and $E=1$  gives case (ii)

$A=D=1$ and $E=0$ gives (vi)

$A=E=1$ and $D=0$ gives (i)

$A=D=1$ and $E=2$ gives (iii)

\n If $B=C=i$ then

$A=1$, $D=-1$ and $E=2$ gives (iv).
\end{D}

By analogy to the K\"ahler case, the reduction of the holonomy of $\na$ to
$SU(n)$ is equivalent to the vanishing of the Ricci form \[\rho(X,Y)=\ft12
\sum_{i=1}^n g(R_{XY}e_i,Je_i),\qq X,Y\in\g\] where $R$ is the curvature tensor
of $\na$. The holonomy reduction of Proposition~\ref{holo} has the advantage of
dealing with a positive definite metric, although the connection is not
torsion-free. However, the situation is just as intriguing for a
pseudo-Riemannian metric and the Levi-Civita connection.

\begin{lem} Suppose that $(J,h)$ is an invariant pseudo-K\"ahler structure on a
nilmanifold $M$ of real dimension $2n$. Then the Ricci tensor of $h$ vanishes.
\end{lem}

\begin{D} The hypothesis means that $M$ admits both an invariant complex 
structure $J$ and a closed 2-form $\Om$ for which $h$ (defined by \eqref{gom}
with $h$ in place of $g$) is a pseudo-Riemannian metric. In the presence of a
compatible complex structure, the equation $d\Om=0$ is sufficient to imply that
$D\Om=0$ and thus $DJ=0$, just as in the familiar (positive-definite) K\"ahler
case. We can therefore apply \eqref{beta} with $d\Om=0$ to deduce that
$D\eta=0$, and the result follows.\end{D}

As an example,

\begin{pro} The Lie algebra $(0,0,0,0,13+42,14+23)$ associated to the Iwasawa 
manifold admits a pseudo-K\"ahler metric $h$ for which $D$ has holonomy in
$SU(2,1)$ and $h$ is Ricci-flat but not flat.\end{pro}

\begin{D} We take \[\Om=e^{16}+e^{25}+e^{34},\qq 
Je^1=e^2,\ Je^3=-e^4,\ Je^5=-e^6.\] It is necessary to check that \[d\Om=
-e^1\we de^6-e^2\we de^5=-e^1\we(e^{14}+e^{23})-e^2\we(e^{13}+e^{42})=0\] and 
that \[ J\cdot\Om=(Je^1)\we(Je^6)+(Je^2)\we(Je^5)+(Je^3)\we(Je^4)=\Om.\]
Setting $\a1=e^1+ie^2,\,\a2=e^3-ie^4,\,\a3=e^5-ie^6$ gives $d\a3=\aa1\2$, so
$J$ is integrable and abelian. Observe that $\hJ$ has negative orientation, so
its connected component in $\cC$ does not contain $\pm J_0$. The resulting
pseudo-metric \eqref{gom} assumes the matrix form
\[(h_{ij})=\left(\!\ba{cccccc} 0&0&0&0&1&0\\0&0&0&0&0&-\!1\\0&0&1&0&0&0\\
0&0&0&1&0&0\\1&0&0&0&0&0\\0&-\!1&0&0&0&0\ea \!\right)\] and has signature 
(4,2). 

A computation reveals that
\[\hbox{$R_{1212}=\sum\limits_{m=1}^6h_{1m}R^m{}_{212}=R^5{}_{212}=-2$,}\]
and so the full Riemann curvature is non-zero.\end{D}

A contrasting situation is provided by the Lie algebra
$(0,0,0,0,0,12+34)$. This carries no invariant symplectic form
\cite{SalCSN}, and does not therefore possess an invariant
pseudo-K\"ahler metric. We suspect that the remaining four algebras in
\eqref{list} do admit pseudo-K\"ahler metrics.

A completely different class of Ricci-flat structures on nilmanifolds, indeed
ones with signature $(2n,2n)$ that are not invariant, are discussed in
\cite{FPPNCY}.

\vs

\n{\small Dipartimento di Matematica, Universit\`a di Torino, Via Carlo
Alberto 10, 10124 Torino, Italy}\\\texttt{annamaria.fino@unito.it}\vs

\n{\small Dipartimento di Matematica, Universit\`a di Pisa, Via Buonarroti 2, 56127 Pisa, Italy}
\\\texttt{parton@dm.unipi.it}\vs

\n{\small Dipartimento di Matematica, Politecnico di Torino, Corso Duca
degli Abruzzi 24, 10129 Torino, Italy}\\\texttt{salamon@calvino.polito.it}

\enddocument